\documentstyle[12pt]{article}
\newtheorem{ppp}{Proposition}[section]
\newtheorem{lll}{Lemma}[section]

\begin{document}

\setcounter{page}{0}

\title{
 Birational automorphisms
\\  of algebraic varieties
\\  with a pencil of cubic surfaces}

\author {Aleksandr V. Pukhlikov \\ Max-Planck-Institut f\"ur Mathematik \\
Gottfried-Claren-Stra{\ss}e 26, D-5300 Bonn 3 \\ and
\\ Institute for Systems Studies \\ Prospekt 60-letya Oktyabrya 9, 117312
Moscow }

\date{October 11, 1996}

\thanks{The author was financially supported by the Russian Fundamental
Research Foundation,
project
 96-01-00820 }

\maketitle

\begin{abstract}
It is proved that on a smooth algebraic variety, fibered into  cubic surfaces
over the projective line and
sufficiently ``twisted'' over the base, there is only one pencil of rational
surfaces -- that is, this very
pencil of cubics. In particular, this variety is non-rational; moreover, it can
not be fibered into
rational curves. The proof is obtained by means of the method of maximal
singularities.
\end{abstract}

\section*{Introduction}

The present paper deals with birational maps of smooth algebraic threefolds,
fibered into cubic surfaces -- that is, Del Pezzo surfaces of the degree 3, --
over a rational curve.
In [1] V.A.Iskovskikh formulated a conjecture of ``birational rigidity'' of the
pencils
of Del Pezzo surfaces of the degrees 1,2 and 3, which are ``sufficiently
twisted''
over the base. For the pencils of surfaces of the degrees 1 and 2 this
conjecture
(to be exact, its slight modification) was proved by means of the method of
maximal singularities in [6]. The present paper is a direct continuation of
[6].
It completes the proof of V.A.Iskovskikh's conjecture in full. The structure of
the
paper is analogous to [6], whereas most of it (Sections 4-6) is devoted to the
study
of infinitely near maximal singularities of linear systems. We are mostly
concerned with
the situation where the techniques developed in [6] do not work: when the
maximal singularity
which is to be excluded is contracted to a point lying on a line in a fiber of
the pencil
of cubic surfaces.
To cope with this situation, we need to develop a new technique, based upon
detailed analysis
of effective 1-cycles on the certain specially chosen blow ups of the original
variety.

In the first section of the paper we describe the class of varieties, the
birational type of
which is to be studied. In the second section we formulate the principal result
-- that is, that
there is exactly one pencil of rational surfaces on the varieties from this
class. Here we also
start to prove it. The third section deals with maximal curves. Here we follow
Yu.I.Manin [4] and
V.A.Iskovskikh [1]. In the rest of the paper, which is actually its principal
part, we study
infinitely near singularities. Note that we do not repeat those proofs which
are completely analogous
to the proofs of the corresponding statements in [6]. We just make reference to
the corresponding
arguments in [6]. All the notations in [6] and in the present paper are
compatible.

Thus the proof of V.A.Iskovskikh's conjecture is now complete. This makes it
possible to
``sum up'' the results which were obtained by means of the method of maximal
singularities
during the twenty five years of its existence, starting from the pioneer paper
of
V.A.Iskovskikh and Yu.I.Manin [2], where G.Fano's ideas were for the first time
realized
on the modern rigorous level. So: in the ``stable'', that is, sufficiently
twisted
over the base, cases  we have the complete description of the birational type
of conic
bundles [8,9] and Del Pezzo fibrations. For certain Fano varieties, including
some singular ones,
analogous results were also obtained, see the survey [3]. There is some hope
that in the nearest
time M.Reid, A.Corti and the author will complete the joint paper in which it
is to be proved that
the weighted Fano hypersurfaces of the index 1 (there are 95 families of them)
are birationally rigid. All this makes it not completely unrealistic to hope
that the method of
maximal singularities can be developed to be powerful enough to ensure the
complete
birational classification  and description of the groups of birational
automorphisms of
varieties with negative Kodaira dimension, at least in dimension three.

The paper was completed at the beginning of my stay at Max-Planck-Institut
f\"ur
Mathematik. I would like to express my gratitude to the staff of the Institute
for their
hospitality.

\section{Varieties with a pencil of cubic surfaces}

We study smooth three-dimensional projective varieties
$V$
over the field of complex numbers
${\bf C}$,
permitting a morphism
$\pi:V\to{\bf P}^1$,
every fiber of which is an irreducible reduced Del Pezzo surface of the degree
3, whereas its generic fiber
$F_{\eta}$
is a smooth Del Pezzo surface of the degree 3 over the non-closed field
${\bf C}({\bf P}^1)$
with the Picard group
$\mathop{\rm Pic} F_{\eta}={\bf Z} K_{F_{\eta}}$.
In particular,
$\mathop{\rm Pic} V={\bf Z} K_V\oplus {\bf Z} F$,
where
$F$
is a fiber of the morphism
$\pi$.
Obviously,
$$
(K^2_V\cdot F)=3.
$$

The relatively ample sheaf
${\cal O}_V(-K_V)$
is generated by its sections on each fiber, and thus determines the inclusion
$$
V\hookrightarrow X,
$$
where
$X={\bf P}({\cal E})$,
${\cal E}={\cal O}\oplus{\cal O}(a_1)\oplus{\cal O}(a_2)\oplus{\cal O}(a_3)$
is a locally free sheaf of rank 4 on
${\bf P}^1$, $0\leq a_1\leq a_2\leq a_3.$
This sheaf coincides, up to twisting by an invertible sheaf, with
$\pi_*{\cal O}_V(-K_V)$,
see [1]. We denote the projection of
 $X$
onto
${\bf P}^1$
by
$\pi$,
too.

Let
$G$
be a fiber of the morphism
$\pi:X\to{\bf P}^1$.
Then
$\mathop{\rm Pic} X={\bf Z} L\oplus{\bf Z} G$,
where
$L$
is the class of the tautological sheaf
${\cal O}_{X/{\bf P}^1}(1)$,
that is,
$\pi_*{\cal O}_X(L)={\cal E}$.
Note that the sheaf
${\cal O}_X(L)$
is generated by global sections. In other words,
$|L|$
is a free linear system.
The class of the variety
$V$
as a divisor on
$X$
is equal to
$3L+mG$.
The canonical class of the variety
$V$
is equal to the restriction onto
$V$
of the class
$$
-L+(m+a_1+a_2+a_3-2)G.
$$

The group of 1-cycles on
$V$
can be written down as
$$
A^2(V)={\bf Z} s\oplus{\bf Z} f,
$$
where
$s$
is the class of a section,
$f$
is the class of a line in a fiber:
$$
(-K_V\cdot f)=(F\cdot s)=1, (F\cdot f)=0.
$$

In the present paper we study those varieties, which satisfy the following
condition:
\newline\newline
the class of 1-cycles
$MK^2_V-f$
is not effective for any
$M\in{\bf Z}$.
\newline\newline

We shall refer to this condition as the
``$K^2$-condition''.

This requirement is satisfied,
if
$V$
is ``sufficiently twisted'' over the base
${\bf P}^1$.
It is satisfied, anyway, if the following inequality holds:

$$
(K^2_V\cdot L)\leq 0,
$$
because the system
$|L|$
is free. Easy computations lead us to the inequality
$$
5m\geq 12-3(a_1+a_2+a_3).
$$
In particular, if
$X={\bf P}^1\times{\bf P}^2$,
i.e.
$a_1=a_2=a_3=0$,
then we get the condition
$m\geq 3$,
that is,
$V$
is given in
$X$
by an equation of the form
$$
p_{3,m}(x_0,x_1,x_2,x_3;u,v)=0
$$
of bidegree
$(3,m), m\geq 3$.

Our
$K^2$-condition
essentially coincides with the
hypothesis of V.A.Is\-kovskikh's conjecture [1],
whereas for the case
$X={\bf P}^1\times {\bf P}^2$
the coincidence is exact. For this reason we identify our
theorem, which is proved below, with the above-mentioned conjecture.

Besides, for some purely technical reasons, we assume that the following
condition
of general position holds:
\newline\newline
if the fiber
$\pi^{-1}(t)\subset V$
over a point
$t\in{\bf P}^1$
is a singular cubic surface, then it has exactly one
singular point, which is a non-degenerate quadratic singularity. Moreover,
exactly six lines on the cubic
$\pi^{-1}(t)=F_t$
pass through this point.\newline\newline

There is no doubt that the theorem proved below is true
whether this assumption is  satisfied or not. However, in the general case
one has to consider a lot of particular geometric configurations and to adjust
the constructions of Sections 4-6 to them. The author has not
completed this job, that is why we restrict ourselves by the situation
of general position.

As it was done in [6], we call an irreducible curve
$C\subset V$
a {\it horizontal} one,
if
$\pi(C)={\bf P}^1$,
and a {\it vertical}
one, if
$\pi(C)$
is a point.
We define a vertical and a horizontal effective
1-cycles, respectively.
The
{\it degree} of
a horizontal 1-cycle
$Z$
is equal to
$(F\cdot Z)$,
of a vertical one -- to
$(-K_V\cdot Z)$.
In both cases it is denoted by
$\mathop{\rm deg} Z$.

The fiber of
$V$
over a point
$t\in{\bf P}^1$
-- a cubic surface -- is denoted by the symbol
$F_t$
or
$F_x$,
if
$x\in F_t$.

\begin{ppp}
Assume
$C\subset V$
to be an irreducible curve,
$x\in C$
to be a point on it. Then
$$
\mathop{\rm mult}\nolimits_xC\leq\mathop{\rm deg} C.
$$
\end{ppp}

{\bf Proof:} this is obvious.

{\bf Definition.}
A line
$L\subset V$
(that is, a vertical curve of the degree 1; in other words, a true line in
${\bf P}^3$)
is said to be
{\it non-special}, if
$L\cap\mathop{\rm Sing} F=\emptyset$,
where
$F\supset L$
is the fiber of the morphism
$\pi$,
containing
$L$.
It is said to be
{\it special},
otherwise. In other words, in the special case the line contains the double
point of the fiber.

\begin{ppp}
Let
$L\subset V$
be a special line,
$F\supset L$
be the corresponding fiber -- a singular cubic surface.
Let
$\sigma:\widetilde V\to V$
be the blowing up of
$L$,
$\widetilde F$
be the proper inverse image of the surface
$F$
on
$\widetilde V$.
Then
$\widetilde F$
is a smooth surface, isomorphic to the blow up of
$F$
at the singular point.
\end{ppp}

{\bf Proof:} straightforward local computations.

Let us also define the ``intersection index'' of two curves on a variety of
arbitrary dimension.

Assume
$R$
to be a nonsingular projective curve on an algebraic variety
$Y$,
$R\cap\mathop{\rm Sing} Y=\emptyset$,
and
$C\subset Y$
to be an arbitrary irreducible curve,
$C\neq R$.

{\bf Definition.}
The
{\it intersection index}
$$
(C\cdot R)
$$
of the curves
$C$
and
$R$
is the integer given by one of the two equivalent constructions:

1)
$(C\cdot R)=\sum\limits_{x\in C\cap R}\mathop{\rm mult}\nolimits_xC$,
where the sum is taken over all the points of intersection of
$C$
and
$R$,
{\it including infinitely near ones};

2)
$(C\cdot R)=(\widetilde C\cdot E)$,
where
$\sigma:\widetilde Y\to Y$
is the blowing up of the curve
$R$
with the exceptional divisor
$E\subset\widetilde Y$,
and
$\widetilde C\subset\widetilde Y$
is the proper inverse image of the curve
$C$.

We define the intersection index
$(Z\cdot R)$
for any 1-cycle
$Z$
on
$Y$,
which does not contain
$R$
as  a component, by linearity.
If
$S\subset Y$
is a smooth projective surface, containing
$R$
and the support of
$Z$,
then we get the intersection index
$(Z\cdot R)$
on
$S$
in the usual sense.

Let us introduce one more notation: the multiplicity of the curve
$C$
in a 1-cycle
$Z$
is denoted simply by
$$
\mathop{\rm mult}\nolimits_C Z.
$$

\section{Formulation of the principal result.\protect\\ Start of the proof}

Now fix a variety
$V$,
for which the
$K^2$-condition
is satisfied. Let
$V'$
be a projective threefold, nonsingular in codimension 1
and fibered into rational surfaces
$$
\pi':V'\to{\bf P}^1,
$$
Let
$F'$
be a fiber of this pencil. Our principal result is the following
\newline\newline
{\bf Theorem.}
{\it Any birational map}
$$
\chi:V-\,-\,\to V'
$$
{\it (provided there are any) is fiber-wise with respect to the pencils}
$\pi,\pi'$,
{\it that is,}
$$
\chi^{-1}(|F'|)=|F|.
$$
{\it In other words, there is an isomorphism of the base}
$$
\alpha:{\bf P}^1\to{\bf P}^1,
$$
{\it such that}
$\pi'\circ\chi=\alpha\circ\pi$.
\newline\newline

{\bf Remark.}
In fact, we could assume
$|F'|$
to be any pencil of surfaces of negative Kodaira dimension.
It is not necessary to make any changes at all in the proof
given below to make it work in this case. But this, formally more general,
statement gives actually no additional information about the
birational type of the variety
$V$.
So we formulate our theorem in the same manner as it was done in [1].

{\bf Corollary.}
(i) $V$
{\it has only one pencil of rational surfaces. In particular, }
$V$
{\it is not rational and has no structures of a conic bundle.}

(ii) {\it The quotient group of the group}
$\mathop{\rm Bir} V$
{\it of birational automorphisms of the variety}
$V$
{\it by the normal subgroup of birational automorphisms,
preserving the fibers of}
$\pi$
{\it (which is isomorphic to the group}
$\mathop{\rm Bir} F_{\eta}$
{\it of birational automorphisms of the generic fiber),}
{\it is finite, generically trivial.}
\newline\newline

The theorem was formulated by V.A.Iskovskikh as a conjecture [1]. The
corollaries
(i)
and
(ii)
follow from the theorem in an obvious way. They were also
formulated in the above-mentioned paper of V.A.Iskovskikh.

\subsection*{Start of the proof}

Let
$|\chi|$
be the proper inverse image of the pencil
$|F'|$
on
$V$.
There exist integers
$n=n(\chi)\in{\bf Z}_+$
and
$l\in{\bf Z}_+$
such that
$$
|\chi|\subset|-n(\chi)K_V+lF|.
$$
If
$n(\chi)=0$,
then
$l=1$
and the theorem is true.

Starting from this moment, we assume that
$n(\chi)\geq 1$.
Obviously, this is equivalent to
$\chi$
being not fiber-wise. We show below that this
assumption leads to a contradiction.

The fact that the pencil
$|\chi|$
has no fixed components implies, in accordance with the
$K^2$-condition, that
$l\geq 0$
(see [6]).

\subsection*{Adjunction break condition}

We use the language of discrete valuations in the form of
[5,7].
The centre of a valuation
$\nu\in{\cal N}(V)$
is denoted by
$Z(V,\nu)$.
If
$Z(V,\nu)$
is a point, then the fiber of the pencil
$|F|$,
containing this point, is denoted by
$F_{\nu}$.
Abusing the notations, we sometimes write
$T$
instead of
$\nu$,
if
$\nu=\nu_T$,
where
$T$
is a prime Weyl divisor on some model
$\widetilde V$
of the field
${\bf C}(V)$,
$T\not\subset\mathop{\rm Sing}\widetilde V$,
realizing the discrete valuation
$\nu$.

For a valuation
$\nu\in{\cal N}(V)$
set
$$
e(\nu)=\nu(|\chi|)-n\delta(\nu),
$$
where
$\delta(\nu)=K(V,\nu)$
is the canonical valuation (discrepancy) of
$\nu, n=n(\chi)$.
The valuations for which
$\nu(|\chi|)>0$
are said to be
{\it singularities} of the system
$|\chi|$.
A discrete valuation for which
$e(\nu)>0$
is said to be a
{\it maximal singularity}.

In the assumptions of the theorem we get

\begin{ppp}
A maximal singularity does exist. Moreover, one of the following two cases
takes place:

(i)  there is a maximal singularity
$\nu\in{\cal N}(V)$
such that its centre
$Z(V,\nu)$
on
$V$
is a curve;

(ii) there is a finite set of maximal singularities
${\cal M}\subset{\cal N}(V)$,
the centres
$Z(V,\nu)=x(\nu)$
of which are points on
$V$,
and, moreover, the following inequality holds
$$
\sum\limits_{t\in{\bf P}^1}
\left(
\max\limits_{\{\nu\in{\cal M}|
x(\nu)\in F_t=F_{\nu}\}}
\frac{e(\nu)}{\nu(F_t)}
\right)>l.
$$
\end{ppp}

The {\bf proof} was given in
[6].
The structure of the rest of the paper is similar to
[6]:
first
(in the next section),
we study the maximal curves
(the case (i) of the Proposition).
The principal part of our proof is concentrated in Sections 4-6, where we show
that the case
(ii)
of the Proposition is impossible -- it leads to a contradiction.

\section{Maximal curves}

Here we follow Yu.I.Manin
[4]
and V.A.Iskovskikh
[1].

If the centre of a maximal singularity
$\nu$
is a curve
$C\subset V$,
then it is easy to see that the curve
$C$
itself is already maximal:
$$
\mathop{\rm mult}\nolimits_C|\chi|>n(\chi).
$$

\subsection*{Case 1: $C$ is horizontal}

First of all, we restrict the linear system
$|\chi|$
onto the generic fiber
$F$
and obtain the inequality
$\mathop{\rm deg} C\leq 2$
(by the arguments, similar to the corresponding ones in [6], Section 4).

Let us start with the case
$\mathop{\rm deg} C=1$,
that is,
$C$
is a section of the morphism
$\pi$.
For a general point
$t\in{\bf P}^1$
take a general line
$L\subset{\bf P}^3=G_t$,
passing through the point
$C\cap F_t$.
This line intersects the cubic surface
$F_t$
at two more different points
$x,y$.
Set
$$
\tau_C(x)=y.
$$
Obviously, by means of this construction the birational involution
$\tau_C\in\mathop{\rm Bir} F_{\eta}\subset\mathop{\rm Bir} V$
is defined. Let
$\alpha:V^*\to V$
be the blowing up of the curve
$C$,
$E=\alpha^{-1}(C)$
be the exceptional divisor,
$\mathop{\rm Pic} V^*={\bf Z} h\oplus{\bf Z} e\oplus{\bf Z} F, h=-K_V$.

\begin{lll}
The birational involution
$\tau_C$
extends to a biregular involution of an invariant open set
$V^*\backslash Y$, $\mathop{\rm codim} Y\geq 2$.
Its action on
$\mathop{\rm Pic} V^*/{\bf Z} F\cong{\bf Z}\bar h\oplus{\bf Z}\bar e$
is given by the relations
$$
\tau^*_C\bar h= 3\bar h-4\bar e,
$$
$$
\tau^*_C\bar e= 2\bar h-3\bar e.
$$
\end{lll}

{\bf Proof.} See
[4].

Now let us
``untwist'' the curve
$C$.
Consider the composition
$\chi\circ\tau_C:V-\,-\,\to V'$.

\begin{lll}
The following relation holds
$$
n(\chi\circ\tau_C)=
3n(\chi)-2\nu_C(\chi)<n(\chi),
$$
Besides, the curve
$C$
is no more maximal for the composition
$\chi\circ\tau_C$.
\end{lll}

{\bf Proof.} The linear system
$|\chi\circ\tau_C|$
is the proper inverse image of
$|\chi|$
with respect to
$\tau_C$.
Respectively, the proper inverse image of the linear system
$|\chi\circ\tau_C|$
on
$V^*$
can be obtained by applying
$\tau_C$
to the proper inverse image of the linear system
$|\chi|$
on
$V^*$.
Thus
$$
n(\chi\circ\tau_C)h-
\nu_C(\chi\circ\tau_C)e=
$$
$$
=\tau^*_C\left(
n(\chi)h-\nu_C(\chi)e
\right)=
\left(3n(\chi)-2\nu_C(\chi)
\right)h+\dots,
$$
and we are done. Q.E.D.
\newline\newline

Now assume that
$\mathop{\rm deg} C=2$,
i.e.
$C$
is a bisection of the morphism
$\pi$.
We define the involution
$\tau_C$
by its action on the generic fiber
$F$
in the following way
(see [4]).
Set
$\{a,b\}=C\cap F$.
Let
$q=L_{ab}\cap F$
be the third point of intersection of the line in
${\bf P}^3$,
joining the points
$a$
and
$b$,
with the cubic surface
$F$.
These points
$q$
sweep out a curve
$C^*\subset V$,
which is a section of the morphism
$\pi$,
i.e.
$q=C^*\cap F$.
The pencil of planes
$P$
in
${\bf P}^3$,
containing the line
$L_{ab}$,
generates the pencil of elliptic curves
$Q_P=P\cap F$
on the surface
$F$.
Set
$$
\tau_C\left|_{Q_P}(x)\right.=y,
$$
where
$$
x+y\sim 2q
$$
on
$Q_P$,
i.e.
$\tau_C$
is the reflection from the point
$x$
on the elliptic curve
$Q_P$.
Thus the involution
$\tau_C\in\mathop{\rm Bir} F_{\eta}\subset\mathop{\rm Bir} V$
is defined.

Let
$\alpha:V^*\to V$
be the blowing up of the curve
$C$,
$E=\alpha^{-1}(C)$
be the exceptional divisor,
$\mathop{\rm Pic} V^*={\bf Z} h\oplus{\bf Z} e\oplus{\bf Z} F, h=-K_V$.

\begin{lll}
The birational involution
$\tau_C$
extends to a biregular involution of an invariant open set
$V^*\backslash Y$, $\mathop{\rm codim} Y\geq 2$.
Its action on
$\mathop{\rm Pic} V^*/{\bf Z} F\cong{\bf Z}\bar h\oplus{\bf Z}\bar e$
is given by the relations
$$
\tau^*_C\bar h= 5\bar h-6\bar e,
$$
$$
\tau^*_C\bar e= 4\bar h-5\bar e.
$$
\end{lll}

{\bf Proof:} straightforward computations.

Now consider the composition
$\chi\circ\tau_C:V-\,-\,\to V'$.

\begin{lll}
The following relation holds
$$
n(\chi\circ\tau_C)=
5n(\chi)-4\nu_C(\chi)<n(\chi).
$$
Moreover, the curve
$C$
is no more maximal for the composition
$\chi\circ\tau_C$.
\end{lll}

{\bf Proof: } straightforward computations, similar to the  above  ones.

The computations which were just performed are none else but the well
known constructions of the two-dimensional birational geometry over non-closed
fields
[4],
translated into the language of a threefold, fibered over
${\bf P}^1$.

\subsection*{Case 2: $C$ is vertical}

It was proved in
[1], that this case does not realize. The most easy way to show it is by means
of the techniques developed in
[5]. First of all, it is easy to see
that
$C\subset F\subset {\bf P}^3$
is a line or a conic.
Furthermore, let
$x\in {\bf P}^3$
be a general point,
$S(x)$
be the cone over
$C$
with the vertex
$x$.
Then
$$
C\cup R(x)=F\cap S(x),
$$
where
$R(x)$
is the residual curve. It was shown in
[5]
that
$(C\cdot R(x))=\mathop{\rm deg} R(x)$
(in the sense of the definition of the ``intersection index'', given in Section
1).
This fact together with the inequality
$\mathop{\rm mult}\nolimits_C|\chi|>n$
implies that
$R(x)$
is a basic curve of the linear system
$|\chi|$.
Consequently,
$F$
is a fixed component of the pencil
$|\chi|$.
Contradiction. Q.E.D.

Summing up all these results, we get

\begin{ppp}
There exists a fiber-wise birational automorphism
$$
\chi^*=\tau_{C_1}\circ\dots\circ\tau_{C_k}
$$
such that the linear system
$|\chi\circ\chi^*|$
has no maximal curves on
$V$.
\end{ppp}

All that we need to show now in order to prove our Theorem, is that the second
case
(case (ii))
of Proposition 2.1
does not realize.

Starting from this moment, we assume thus that this very case takes place.

Let us show that this assumption leads to a contradiction.

\section{Infinitely near singularities I.\protect\\ Existence of a line }

\subsection*{The supermaximal singularity}

The symbols
$D_1,D_2\in|\chi|$
stand for generic divisors of our linear system. Consider the effective 1-cycle
$$
Z=(D_1\bullet D_2).
$$
It can be decomposed into the vertical and horizontal components:
$$
Z=Z^v+Z^h,
$$
whereas for the vertical cycle
$Z^v$
we get the decomposition
$$
Z^v=
\sum_{t\in{\bf P}^1}Z^v_t,
$$
where the support of the cycle
$Z^v_t$
lies in the fiber
$F_t$.

\begin{ppp}
There is a maximal singularity
$\nu=\nu_T\in {\cal M}$
such that
$x=x(\nu_T)=Z(V,\nu)$
is a point in a fiber
$F=F_t$,
lying on a line
$L\subset F$.
Moreover, if
$$
Z^v_t=C+kL,
$$
where the effective 1-cycle
$C\subset F$
does not contain the line
$L$,
then
$$
(C\cdot L)<
\frac{4ne}{\nu_T(F)}.
$$
\end{ppp}

{\bf Proof.}
It was proved in
[6]
that the fact that the linear system
$|\chi|$
has no maximal curves implies the existence of a {\it supermaximal singularity}
$\nu=\nu_T\in{\cal M}$,
which satisfies the inequality
$$
\mathop{\rm deg} Z^v_t<
\frac{2dne(T)}{\nu(F_t)},
$$
where
$x(\nu)\in F_t$
and
$d$
is the degree of the generic fiber, that is, the Del Pezzo surface
$F_{\eta}$.
In our case
$d=3$.
Let us show that the supermaximal singularity
$\nu$
satisfies our proposition.

As it was done in
[6],
we write
$x$
instead of
$x(\nu_T)$,
$e$
instead of
$e(T)$, $F$
instead of
$F_t$,
$Z^v$
instead of
$Z^v_t$.

\begin{lll}
The following inequality is true:
$$
\mathop{\rm mult}\nolimits_xZ^v\geq
\frac{4ne}{\nu_T(F)}.
$$
\end{lll}

{\bf Proof.}
Assume the converse. Now, repeating the arguments of Section 5 in [6] word for
word,
we come to a contradiction, for the only fact, upon which they were based, was
exactly
the inequality

$$
\mathop{\rm mult}\nolimits_xZ^v<
\frac{4ne}{\nu_T(F)}.
$$
If it is true, then the supermaximal singularity just cannot exist.
Q.E.D. for the Lemma.

\subsection*{Existence of a line}

\begin{lll}
There is at least one line
$L\subset F\subset {\bf P}^3$
passing through the point
$x\in F$.
\end{lll}

{\bf Proof.}
Assume the converse. Then the point
$x$
is a smooth point of the cubic surface
$F\subset{\bf P}^3$.
Moreover, the curve
$R=T_xF\cap F$
is irreducible, its degree is equal to 3 and
its multiplicity at the point
$x$
is equal to 2 exactly. If
$C\subset F$
is any other curve, then
$$
\mathop{\rm deg} C=(C\cdot R)\geq
(C\cdot R)_x\geq
2\mathop{\rm mult}\nolimits_xC.
$$
Thus for any curve
$Q\subset F$
we get the inequality
$$
\mathop{\rm mult}\nolimits_xQ\leq
\frac23\mathop{\rm deg} Q.
$$
Consequently,
$$
\mathop{\rm mult}\nolimits_xZ^v\leq
\frac23\mathop{\rm deg} Z^v<
\frac{4ne}{\nu(F)}.
$$
Contradiction. Q.E.D.

Now let
$x\in F$
be a smooth point. There exist
$k$
lines, lying on
$F$,
$1\leq k\leq 3$
and passing through
$x$. If
$\mathop{\rm mult}\nolimits_x(F\cap T_xF)=2$,
then
$k\leq 2$
and for any curve
$C\subset F$,
which is different from these
$k$
lines, we get the inequality
$$
2\mathop{\rm mult}\nolimits_x C\leq \mathop{\rm deg} C.
$$
If
$\mathop{\rm mult}\nolimits_x(F\cap T_xF)=3$,
then this inequality can be strengthened:
$$
3\mathop{\rm mult}\nolimits_xC\leq \mathop{\rm deg} C.
$$

\subsection*{The case of a single line}

Assume that there is only one line on
$F$,
passing through
$x$.
In this case the point
$x$
is smooth on
$F$
and
$T_xF\cap F=L+Q$,
where
$Q\subset F$
is a smooth conic. The arguments given above show that for any curve
$C\subset F$,
$C\neq L$,
the following inequality takes place:
$$
\mathop{\rm mult}\nolimits_xC\leq
\frac12\mathop{\rm deg} C.
$$
Write down
$Z^v=C+kL$,
where
$C$
is an effective 1-cycle, not containing
$L$.
Now
$$
k+\frac12\mathop{\rm deg} C\geq
\frac{4ne}{\nu(F)},
$$
$$
k+\mathop{\rm deg} C<
\frac{6ne}{\nu(F)}.
$$
This implies that
$$
\mathop{\rm deg} C<
\frac{4ne}{\nu(F)}.
$$
Since
$(C\cdot L)\leq \mathop{\rm deg} C$,
Proposition 4.1 is proved in this case.

\subsection*{The case of two lines}

Assume that there are exactly two lines,
$L_1$
and
$L_2$
on
$F$,
passing through
$x$.
In this case
$x$
is a smooth point of
$F$
and
$T_xF\cap F=L_1+L_2+L_3$,
where
$L_3$
is a line, different from
$L_1$, $L_2$,
$x\not\in L_3$.
Note that one of the points
$L_i\cap L_3$,
$i=1,2$
can be singular on
$F$.
Assume at first that this is not the case. Write down
$Z^v=Q+k_1L_1+k_2L_2+k_3L_3$
and set
$d=\mathop{\rm deg} Q$,
$d_i=(Q\cdot L_i)$,
$m=\mathop{\rm mult}\nolimits_xQ$.
We get the following four inequalities
$$
2m\leq d,
$$
$$
k_1+k_2+k_3+d<
\frac{6ne}{\nu(F)},
$$
$$
d_1+d_2+d_3=d,
$$
$$
k_1+k_2+m\geq
\frac{4ne}{\nu(F)}.
$$
It is easy to see that there is
$i\in\{1,2\}$
such that for
$\{j\}=\{1,2\}\backslash\{i\}$
$$
k_j+k_3+d_i<
\frac{4ne}{\nu(F)}.
$$
Indeed, otherwise for any
$i\in\{1,2\}$
we have the opposite inequalities.
Put them up together and add the fourth inequality. We get
$$
2(k_1+k_2+k_3)+m+d_1+d_2\geq
\frac{12ne}{\nu(F)}.
$$
This contradicts  the second inequality. So we may assume that
$\nu(F)(k_2+k_3+d_1)<4ne.$
Setting
$C=Q+k_2L_2+k_3L_3$,
we get exactly our proposition.

Finally, if the point
$L_1\cap L_3$
is singular on
$F$,
then, instead of the equality
$d_1+d_2+d_3=d$
one should use the inequality
$d_1+d_2\leq d$,
where
$d_1=(Q\cdot L_1)$
is  taken in the sense of Section 1 of the present paper. It is this
very inequality upon which our arguments are actually based.

\subsection*{The case of three lines}

Assume that there are exactly three lines
$L_i$,
$i=1,2,3$,
on
$F$,
passing through
$x$.
Again
$x$
is smooth on
$F$
and
$T_xF\cap F=L_1+L_2+L_3$.
Write out the 1-cycle
$Z^v=Q+k_1L_1+k_2L_2+k_3L_3$
and set
$d=\mathop{\rm deg} Q$,
$d_i=(Q\cdot L_i)$,
$m=\mathop{\rm mult}\nolimits_xQ$.
Again we get a set of inequalities,
$$
3m\leq d,
$$
$$
k_1+k_2+k_3+d<
\frac{6ne}{\nu(F)},
$$
$$
d_1+d_2+d_3=d,
$$
$$
k_1+k_2+k_3+m\geq
\frac{4ne}{\nu(F)}.
$$
By means of the same elementary arithmetic as above we get that
for some
$i\in\{1,2,3\}$
and
$\{j,l\}=\{1,2,3\}\backslash\{i\}$
the following inequality is true:
$$
k_j+k_l+d_i<
\frac{4ne}{\nu(F)}.
$$
Setting
$L=L_i$,
$C=k_jL_j+k_lL_l+Q$,
we get our proposition.

\subsection*{The case of six lines}

In this case
$x\in F$
is a double point
(a non-degenerate elementary singularity).
Here it is not enough just to make reference to Lemma 4.1.
It is necessary to retrace the arguments of [6]
(Section 5) in details.

Recall our principal notations:
$$
\begin{array}{cccc}
\displaystyle
\varphi_{i,i-1}: & V_i & \to & V_{i-1} \\
\displaystyle
            & \bigcup & & \bigcup \\
\displaystyle
            & E_i & \to & B_{i-1}
\end{array}
$$
is the resolution of
$\nu$
[5,7], that is,
$\varphi_{i,i-1}$
blows up the irreducible cycle
$B_{i-1}$
-- the centre
$Z(V_{i-1},\nu)$
of the valuation
$\nu$
on
$V_{i-1}$,
$E_i=\varphi^{-1}_{i,i-1}(B_{i-1})$
is the exceptional divisor,
$1\leq i\leq K$,
$\nu_{E_K}=\nu$.
The first
$L\leq K$
centres of the blowing ups are points,
after that we blow up curves covering one another.
We equip the set of indices
$\{1,\dots,K\}$
with the natural oriented graph structure
[2,3,5-7],
$p(i.j)$
stands for the number of paths from
$i$
to
$j$,
if
$i\neq j$,
and
$p(i,i)=1$.
Among the six lines
$L_i\ni x$
we choose one
(let it be
$L_1$),
such that
$$
B_1\not\in L^1_i
$$
for
$i\neq 1$,
if
$B_1$
is a point. Such a line does exist. Now write down explicitly
$$
Z^v=kL_1+R+Q,
$$
where
$R$
consists of the multiple lines
$L_i$,
$i\neq 1$,
whereas the 1-cycle
$Q$
does not contain any line, passing through
$x$.
We define the integers
$M$
and
$N$,
by requiring that
$$
B_{i-1}\in L^{i-1}_1, i=1,\dots,M,
$$
$$
B_{i-1}\in F^{i-1}, i=1,\dots,N,
$$
$M\leq N\leq L$,
and set
$$
q_i=
\mathop{\rm mult}\nolimits_{B_{i-1}}Q^{i-1}
$$
for
$i=1,\dots,N$.
Obviously,
$q_1\geq q_2\geq\dots $
and
$$
d=\mathop{\rm deg} Q\geq q_1+q_2+\dots+q_N
$$
(since
$Q$
does not contain the line
$L_1$
as a component).
Now, if the following inequality holds:
$$
\left(
\sum^M_{i=1}p_i
\right)
k+
p_1\mathop{\rm deg} R+
\sum^N_{i=1}p_iq_i<
4ne,
$$
then we immediately get a contradiction by means of the
Iskovskikh-Manin's techniques, repeating the arguments of [6],
Section 5 word for word. Consequently, the opposite inequality holds.
Since
$2q_2\leq d$,
then the following inequality is true:
$$
\left(
\sum^M_{i=1}p_i
\right)
k+
\frac12\left(2p_1+
\sum^N_{i=2}p_i\right)
(\mathop{\rm deg} R+\mathop{\rm deg} Q)\geq
4ne.
$$
Comparing it with the inequality
$$
\left(
2p_1+
\sum^N_{i=2}p_i\right)
(k+\mathop{\rm deg} R+\mathop{\rm deg} Q)<
6ne,
$$
which is true by definition of a supermaximal singularity, we get the following
estimate:
$$
k>
\frac{2ne}{\nu(F)}.
$$
Setting
$C=R+Q$,
we complete the proof of the proposition.

Note that some more delicate arguments
(which are based upon the properties of the integers
$p_j$ only)
make it possible to obtain a stronger estimate of the integer
$(C\cdot L)$
in the last case. However, to prove our main result, we need just the
inequality of
Proposition 4.1.

\section{Infinitely near maximal singularities II. \protect\\ The basic
construction}

\subsection*{Description of the basic construction}

An infinite series of blow ups
$$
\begin{array}{cccc}
\displaystyle
\sigma_i: & V^{(i)} & \to & V^{(i-1)} \\
\displaystyle
            & \bigcup & & \bigcup \\
\displaystyle
            & E^{(i)} & \to & L_{i-1},
\end{array}
$$
$i\geq 1$,
starting from
$V^{(0)}=V$,
where
$L_{i-1}$
is the centre of the
$i$-th
blow up, and
$E^{(i)}=\sigma^{-1}_i(L_{i-1})$
is its exceptional divisor,
$L_0=L$,
is said to be a {\it staircase, associated with the line}
$L$,
or, simply, an $L$-{\it staircase}, if the following conditions are
satisfied:\newline\newline
$L_i$
is a curve for all
$i\in{\bf Z}_+$,
$E^{(i)}$
is a ruled surface of the type
${\bf F}_1$
over
$L_{i-1}$
and
$L_i\subset E^{(i)}$
is the exceptional section
(i.e. the
(-1)-curve).\newline\newline

Obviously, by this definition the staircase is unique. Just below we show that
it exists.
Its segment, consisting of the blow ups
$\sigma_i$
for
$1\leq i\leq M$,
is said to be a (finite) staircase of the length
$M$.

It is convenient to prove the existence of the staircase together with some of
its properties.

For conveniency of notations set
$E^{(0)}$
to be the fiber
$F$
of the morphism
$\pi$,
which contains
$L$.
The operation of taking the proper inverse image on the
$i$-th step
(i.e. on
$V^{(i)}$)
is denoted by adding the bracketed upper index
$i$.
For instance, the proper inverse image of the surface
$E^{(i)}$
on
$V^{(j)}$
for
$j\geq i$
is written down as
$E^{(i,j)}$.
Set also:

$s_i$
to be the class of
$L_i$
in
$A^2(V^{(i)})$,
$s_0=f$;

$f_i\in A^2(V^{(i)})$
to be the class of the fiber of the ruled surface
$E^{(i)}$
over a point
$\in L_{i-1}$.

Abusing our notations, we sometimes treat
$s_i$
and
$f_i$
as numerical classes of curves on the ruled surface
$E^{(i)}$:
$$
A^1E^{(i)}=
\mathop{\rm Pic} E^{(i)}=
{\bf Z} s_i\oplus{\bf Z} f_i,
$$
so that, in particular, the formulas like
$$
(s_i\cdot s_i)=-1,
$$
$$
(s_i\cdot f_i)=1
$$
make sense.

In these notations we have the following

\begin{ppp}
(i) For
$i\geq 2$
the effective 1-cycle
$(E^{(i-1,i)}\bullet E^{(i)})$
is just the irreducible curve
$E^{(i-1,i)}\cap E^{(i)}$.
Its numerical class is equal to
$(s_i+f_i)$.
In particular, this curve does not intersect
$L_i\sim s_i$.
If the line
$L$
is non-special, then this statement is true for
$i=1$,
too. If, on the contrary,
$L$
is special, then the 1-cycle
$(F^{(1)}\bullet E^{(1)})$
is a reducible curve. More exactly, it is the sum of the exceptional section
$L_1$
and the fiber over the singular point of the surface
$F$.

(ii) The following equalities hold:
$$
(E^{(i)})^3=1,
$$
$$
(E^{(i)}\cdot L_i)=0.
$$
Taking into account the isomorphism
$L_i\cong{\bf P}^1$,
we can write down
$$
{\cal N}_{L_i/V^{(i)}}\cong
{\cal O}_{L_i}\oplus{\cal O}_{L_i}(-1).
$$
In this representation the first component is uniquely determined. It
corresponds
to the exceptional section
$L_{i+1}\subset E^{(i+1)}=
{\bf P}({\cal N}_{L_i/V^{(i)}})$.
For the second component we  can take the one-dimensional subbundle,
corresponding exactly to the curve
$E^{(i,i+1)}\cap E^{(i+1)}$.

(iii) The classes
$s_i$
and
$f_i$
satisfy the relations
$$
\sigma^*s_{i-1}=s_i,
$$
$$
\sigma_* f_i=0
$$
for
$i\geq 1$.
\end{ppp}

\subsection*{Proof in the non-special case}
Assuming that
$L\cap\mathop{\rm Sing} F=\emptyset$,
we prove simultaneously the existence of the staircase and Proposition 5.1.

Let us consider the first step of the staircase, that is, the morphism
$$
\sigma_1:
V^{(1)}\to
V^{(0)}=V,
$$
blowing up the line
$L_0=L\subset F$.
We get the exact sequence
$$
0\to
{\cal N}_{L/F}\to
{\cal N}_{L/V}\to
{\cal O}_V(F)\left|_L\right.\to 0,
$$
which can be rewritten down in the following way:
$$
0\to
{\cal O}_L(-1)\to
{\cal N}_{L/V}\to
{\cal O}_L\to 0.
$$

Consequently,
$E^{(1)}$
is a ruled surface of the type
${\bf F}_1$,
$(E^{(1)})^3=1$,
whence
$(E^{(1)}\cdot E^{(1)})\sim (-s_1-f_1)$
and
$(E^{(0,1)}\cdot L_1)$
$=((F-E^{(1)})\cdot s_1)=0$.
Thus all the requirements
(i)-(iii)
of the definition of the staircase are satisfied for the first blow up.

We proceed by induction on
$i\geq 1$.
We get the exact sequence
$$
0\to
{\cal N}_{L_i/E^{(i)}}\to
{\cal N}_{L_i/V^{(i)}}\to
\left.{\cal O}_{V^{(i)}}(E^{(i)})\right|_{L_i}\to 0.
$$
Taking into account the facts which were already proved, this sequence can be
rewritten
down as follows:
$$
0\to
{\cal O}_{L_i}(-1)\to
{\cal N}_{L_i/V^{(i)}}\to
{\cal O}_{L_i}\to
0.
$$
Again this implies that
$E^{(i+1)}={\bf P}({\cal N}_{L_i/V^{(i)}})$
is a ruled surface of the type
${\bf F}_1$
and
$(E^{(i+1)})^3=1$,
so that
$$
E^{(i+1)}\left|_{E^{(i+1)}}\right.\sim
(-s_{i+1}-f_{i+1}).
$$
Thus
$(E^{(i+1)}\cdot L_{i+1})=0$,
(i)
and
(iii)
are satisfied in an obvious way. The proof is complete. Q.E.D.

\subsection*{Proof in the special case}

Assume that the line
$L$
contains the double point
$p\in F$
of the fiber. Again consider the first blow up:
$$
\sigma_1:V^{(1)}\to V^{(0)}=V.
$$
As we have mentioned, the proper inverse image of the fiber
$F^{(1)}=E^{(0,1)}$
is already a non-singular surface. It is easy to see that
$(F^{(1)}\bullet E^{(1)})=L_1+R$,
where
$R=\sigma^{-1}_1(p)$
is the fiber of
$\sigma_1$
over the singular point, whereas
$L_1$
is a certain section of the ruled surface
$E^{(1)}$.
Since
$K_{F^{(1)}}=\sigma^*_1 K_F$,
we get
$$
(L_1\cdot L_1)_{F^{(1)}}=-1,
$$
so that
$(L_1\cdot E^{(1)})=0$
and
$(L_1\cdot F^{(1)})=0$.
Since
$E^{(1)}\left|_{E^{(1)}}=
-F^{(1)}\right|_{E^{(1)}}=
-(L_1+R)$,
we get
$$
(L_1\cdot L_1)_{E^{(1)}}=-1,
$$
so that
$E^{(1)}$
is a ruled surface of the type
${\bf F}_1$
and the conditions
(i)-(iii)
of Proposition 5.1 are satisfied.

The rest of the arguments
(for
$i\geq 2$)
just repeat word for word the non-special case.

The proof of the existence of the staircase and of Proposition 5.1 is complete.
Q.E.D.

{\bf Remarks.}
(i) Since
$E^{(i-1,i)}$
does not intersect
$L_i$
(for
$i\geq 1$
in the non-special and for
$i\geq 2$
in the special case), we get
$$
E^{(i-1,i)}=
E^{(i-1,i+1)}=
\dots
=E^{(i-1,j)}=
\dots
$$
for any
$j\geq i$.
In particular, if
$C\subset E^{(i-1)}$
is a curve, which is not the exceptional section
$L_{i-1}$,
then its proper inverse images on all the varieties
$V^{(j)}$,
$j\geq i$,
are the same:
$$
C^{(i)}=
C^{(i+1)}=
\dots=
C^{(j)}.
$$

(ii) Abusing our notations, we call an irreducible curve
$C\subset E^{(i)}$,
$i\geq 1$,
a {\it horizontal } one, if
$\sigma_i(C)=L_{i-1}$,
and a {\it vertical} one, if
$\sigma_i(C)$
is a point on
$L_{i-1}$.
Respectively, we define horizontal and vertical 1-cycles with the support in
$E^{(i)}$.
The {\it degree} of a horizontal curve
$C$
is equal to
$\mathop{\rm deg} C=$ $\mathop{\rm deg}\sigma_i\left|_C\right.$
$=(C\cdot f_i)$,
the {\it degree } of a vertical curve
$C$
is equal to
$\mathop{\rm deg} C=(C\cdot L_i)=1$.
We define the degree of a horizontal and a vertical 1-cycle with the support in
$E^{(i)}$
as its intersection with
$f_i$
and
$L_i$,
respectively. In particular, the degree of a vertical 1-cycle is just the
number of lines (fibers) in it.
Note that if an effective horizontal 1-cycle
$C$
does not contain the exceptional section
$L_i$
as a component, then its class in
$A^1(E^{(i)})$
or
$A^2(V^{(i)})$
is equal to
$\alpha s_i+\beta f_i$,
where
$\alpha \geq 1$
and
$\beta\geq\alpha$.

(iii) Obviously, the graph of the sequence of the blow ups
$\sigma_i$
is a chain. In particular,
$$
K_{V^{(M)}}=
\sigma^*_{M,0}K_V+
\sum^M_{i=1}
\sigma^*_{M,i}E^{(i)}
$$
(where
$\sigma_{i,j}$,
as always, stands for the composition
$\sigma_{j+1}\circ\dots\circ\sigma_i$)
and the canonical multiplicity of the valuation
$\nu_{E^{(i)}}$
is equal to
$i$.
In the non-special case
$$
\sigma^*_{M,0}F=
F^{(M)}+
\sum^M_{i=1}E^{(i,M)},
$$
whereas in the special case for
$M\geq 2$

$$
\sigma^*_{M,0}F=
F^{(M)}+
E^{(1,M)}+
2\sum^M_{i=2}E^{(i,M)}.
$$

In the special case
$F^{(1)}\cap E^{(1)}$
is equal to the reducible curve
$L_1+R$,
$R=\sigma^{-1}_1(p)$
is the fiber over the double point of the surface
$F$.
There are five more lines on
$F$
besides
$L$
(in accordance with the condition of general position),
passing through
$p$.
Let
$Q$
be one of them,
$Q^{(1)}\subset F^{(1)}$
be its proper inverse image on
$V^{(1)}$.
It is easy to see that the point
$$
\left(
Q^{(1)}\cap E^{(1)}
\right)
\in R
$$
does not lie on the exceptional section
$L_1$.

\section{Infinitely near maximal singularities III.\protect\\ Completing the
proof}

The present section is the key part of the paper. Here we exclude the
supermaximal
singularity. In accordance with what was proved in Section 4, the centre
$x=Z(V,\nu)$
lies on a line
$L\subset F$
and
$$
(C\cdot L)<
\frac{4ne}{\nu(F)},
$$
$$
Z^v=C+kL.
$$

\begin{ppp}
There exists a finite
$L$-staircase
of the length
$M\geq 1$
satisfying the following conditions:

(i) for
$i=0,\dots,M-1$
the centre
$Z(V^{(i)},\nu)$
of the valuation
$\nu$
on
$V^{(i)}$
is a point
$x_i\in L_i$,
$x_0=x$,

(ii) the centre
$Z(V^{(M)},\nu)$
is either:

A) a point
$x_M\not\in L_M, x_M\not\in E^{(M-1,M)}$;

B) the line
$B=\sigma^{-1}_M(x_{M-1})$,
that is, a fiber of the ruled surface
$E^{(M)}$;

C) the point
$x_M=
E^{(M-1,M)}\cap \sigma^{-1}_M(x_{M-1})$.
\end{ppp}

{\bf Proof.}
If
$Z(V^{(i)},\nu)\subset E^{(i)}$,
then
$i=K(V,E^{(i)})\leq K(V,\nu)$.
Consequently, there exists an integer
$M\geq 1$,
such that for
$i=0,\dots,(M-1)$
the condition
(i)
is satisfied, whereas for
$i=M$
it is not satisfied. Now the centre
$Z(V^{(M)},\nu)$
is either a curve
(and in this situation we get the case B)),
or a point
$x_M$,
not lying on
$L_M$
(one of the cases A) and
C)).
Q.E.D. for the Proposition.

Let us fix the just constructed staircase of the length
$M$.

Let us introduce a new parameter
$d_F$,
setting it to be equal to 1, if the case is non-special, and 2, otherwise.
Denote
$E^{(M)}\subset V^{(M)}$
by
$E$
in the cases
A)
and
B),
$E^{(M)}$
by
$E_+$
and
$E^{(M-1,M)}$
by
$E_-$
in the case
C).

\subsection*{Noether-Fano inequality in terms of the staircase}

Denote by
$|\chi|^{(i)}$
the proper inverse image of the system
$|\chi|$
on
$V^{(i)}$
and set
$$
\lambda_i=
\mathop{\rm mult}\nolimits_{L_{i-1}}|\chi|^{(i-1)},
$$
$n\geq\lambda_1\geq\dots$.
Let
$$
\begin{array}{cccc}
\displaystyle
\varphi_{i,i-1}: & V_i & \to & V_{i-1} \\
\displaystyle
            & \bigcup & & \bigcup \\
\displaystyle
            & E_i & \to & B_{i-1},
\end{array}
$$
$i=1,\dots,K$, $V_0=V^{(M)}$,
be the resolution of the valuation
$\nu\in{\cal N}(V^{(M)})$.

Let us introduce some new notations:
$$
\nu_i=\mathop{\rm mult}\nolimits_{B_{i-1}}|\chi|^{i-1}
$$
is the multiplicity of the proper inverse image of the system
$|\chi|$
on
$V_{i-1}$
along the cycle which is to be blown up;

$p_i=p(E_K,E_i)$
is the number of paths in the oriented graph of the valuation
$\nu=\nu_{E_K}$,
leading from
$E_K$
to
$E_i$
(here
$\nu$
is considered as a discrete valuation on the variety
$V_0=V^{(M)}$!);

$N^*=\max\{i|1\leq i\leq K, B_{i-1}\subset E^{i-1}\}$
in the cases
A)
and
B);

$L=\max \{i|1\leq i\leq K, B_{i-1}  \mbox{ is a point } \}$
(so that for
$j\leq L$
$B_{j-1}$
is a point, whereas for
$j\geq L+1$
$B_{j-1}$
is a curve, see
[6,7])
in the cases
A)
and
C);

$N=\min\{N^*,L\}$
in the case
A);

$N=N^*$
in the case
B);

$N^*_{\pm}=\max\{i|1\leq i\leq K, B_{i-1}\subset E^{i-1}_{\pm}\}$
in the case
C),
where the signs
$+$
or
$-$
are chosen to be the same in the right-hand and in the left-hand parts;

$N_{\pm}=\min\{N^*_{\pm},L\}$
(in the right-hand part there is the minimum of two integers);

$\Sigma_0=\sum\limits^L_{i=1}p_i$,
$\Sigma_1=\sum\limits^K_{i=L+1}p_i$
in the cases
A)
and
C);

$\Sigma=\sum\limits^K_{i=1}p_i$
in the case
B);

$\Sigma^*=\sum\limits^{N^*}_{i=1}p_i$,
$\Sigma_*=\sum\limits^N_{i=1}p_i$
in the cases
A)
and
B);

$\Sigma^*_{\pm}=\sum\limits^{N^*_{\pm}}_{i=1}p_i$,
$\Sigma_{\pm}=\sum\limits^{N_{\pm}}_{i=1}p_i$
in the case
C).

Obviously, in these notations we get
$$
\nu(E)=\varepsilon=\Sigma^*,
$$
$$
\nu(E_{\pm})=\varepsilon_{\pm}=\Sigma^*_{\pm}.
$$

In the non-special case
A)
or
B)
for
$M\geq 2$
and in the special case
A)
or
B)
for
$M\geq 3$

$$
\nu(F)=d_F\varepsilon,
$$
whereas in the case
C)
(under the same restrictions on
$M$)

$$
\nu(F)=d_F(\varepsilon_++\varepsilon_-).
$$

In either of the cases
A)
or
B)
we get
$$
\nu(|\chi|)=
\nu_E(|\chi|)\nu(E)+
\nu(|\chi|^{(M)})
$$
and

$$
K(V,\nu)=
K(V,\nu_E)\nu(E)+
K(V^{(M)},\nu),
$$
so that the Noether-Fano inequality takes the form
$$
\sum^K_{i=1}p_i\nu_i=
\varepsilon\sum^M_{i=1}(n-\lambda_i)+
n\sum^K_{i=1}p_i\delta_i+e.
$$
In a similar way, in the case
C)
the Noether-Fano inequality takes the form
$$
\sum^K_{i=1}p_i\nu_i=
\varepsilon_+\sum^M_{i=1}(n-\lambda_i)+
\varepsilon_-\sum^{M-1}_{i=1}(n-\lambda_i)+
n\sum^K_{i=1}p_i\delta_i+e.
$$

\subsection*{Iskovskikh-Manin's techniques}

As always, let
$D^{(M)}_i$,
$i=1,2$,
be the proper inverse images of general divisors from the pencil
$|\chi|$.
Let
$$
Z^{(M)}=
(D^{(M)}_1\bullet D^{(M)}_2)
$$
be the effective 1-cycle of their scheme-theoretic intersection. Set
$$
m_i=
\mathop{\rm mult}\nolimits_{B_{i-1}}(Z^{(M)})^{i-1}
$$
for
$i\leq L$
in the cases
A)
and
C).
In accordance with the Iskovskikh-Manin's techniques
[6,7],
we obtain the following estimate for the case
A):
$$
\sum^L_{i=1}p_im_i\geq
\frac{
(2\Sigma_0n+\Sigma_1n+
\varepsilon\sum\limits^M_{i=1}(n-\lambda_i)+
e)^2}{
\Sigma_0+\Sigma_1}.
$$

For the case
C)
we get the estimate
$$
\sum^L_{i=1}p_im_i\geq
\frac{1}{\Sigma_0+\Sigma_1}\times
$$
$$
\times
\left(
2\Sigma_0n+\Sigma_1n+
\varepsilon_+\sum\limits^M_{i=1}(n-\lambda_i)+
\varepsilon_-\sum\limits^{M-1}_{i=1}(n-\lambda_i)+
e\right)^2.
$$
In the case
B)
we can obviously assert that the line
$B$
comes into the 1-cycle
$Z^{(M)}$
with the multiplicity at least
$$
\sum^K_{i=1}\nu^2_i,
$$
whereas the multiplicities
$\nu_i$
satisfy the inequalities
$$
\nu_i\geq
\sum_{j\to i}\nu_j
$$
(here the resolution of the maximal singularity
$\nu$
is just a sequence of blowing ups of curves, covering each other).
Computing the minimum of this quadratic form under the restrictions specified
above and
taking into account the Noether-Fano inequality, we get
$$
\mathop{\rm mult}\nolimits_BZ^{(M)}\geq
\frac{
\left(
\Sigma n+\varepsilon\sum\limits^M_{i=1}(n-\lambda_i)+e\right)^2}{
\sum\limits^K_{i=1}p^2_i}.
$$

\subsection*{The cycle  $Z^{(M)}$  in terms of the staircase}

Now to complete the proof of our theorem we must get some estimates of the
upper
bounds of the left-hand parts of the three principal inequalities, which were
obtained above.
The computations to be performed are rather tiresome. However, they are quite
clear geometrically.
Coming back to our basic construction -- that is, the staircase,-- let us
introduce some new
terminology and notations, connected with the linear system
$|\chi|$.
First of all, set
$$
z_i=(D^{(i)}_1\cdot D^{(i)}_2)\in A^2V^{(i)}
$$
to be the class of the effective 1-cycle
$$
Z^{(i)}=(D^{(i)}_1\bullet D^{(i)}_2).
$$
On the ``zeroth'' step of our staircase we have the decomposition
$$
Z=Z^v+Z^h.
$$
Let us trace down the changes which the 1-cycle
$Z^{(k)}$
undergoes when
$k$
comes from
$i-1$
to
$i$.
Naturally, instead of the components of the cycle
$Z^{(i-1)}$,
which are different from
$L_{i-1}$,
their proper inverse images come into the cycle
$Z^{(i)}$.
Instead of the curve
$L_{i-1}$,
which is present in
$Z^{(i-1)}$
with some multiplicity
$k_{i-1}$,
the cycle
$Z^{(i)}$
contains an effective sub-cycle with the support in the exceptional divisor
$E^{(i)}$.
Let us break this sub-cycle into three parts:

1)
$C^{(i)}_h$
includes all the curves, which are horizontal with respect to the morphism
$\sigma_i:E^{(i)}\to L_{i-1}$,
and different from the exceptional section
$L_i$,

2)
$C^{(i)}_v$
includes all the vertical curves, that is, the fibers of
$\sigma_i$
over points of the curve
$L_{i-1}$,

3) the exceptional section
$L_i$
with a certain multiplicity
$k_i\in{\bf Z}_+$.

To make our notations look uniform set
$C^{(0)}_h$
to be the part of the cycle
$Z^v$,
which includes all the curves different from
$L$.
Set also
$$
d^{(i)}_{h,v}=\mathop{\rm deg} C^{(i)}_{h,v}
$$
(see Remark
(ii)
in the previous section).
Now we get the following representation of the cycles
$Z^{(i)}$:
$$
Z^{(0)}=Z^h+Z^v=Z^h+C^{(0)}_h+k_0L,
$$
$$
Z^{(1)}=
(Z^h)^{(1)}+
C^{(0,1)}_h+
C^{(1)}_h+
C^{(1)}_v+
k_1L_1,
$$
$$
\dots,
$$
$$
Z^{(i)}=
(Z^h)^{(i)}+
C^{(0,a)}_h+
C^{(1,2)}_h+
C^{(1,2)}_v+
\dots +
$$
$$
+
C^{(i-1,i)}_h+
C^{(i-1,i)}_v+
C^{(i)}_h+
C^{(i)}_v+
k_iL_i.
$$
Here
$a=1$
in the non-special and
$a=2$
in the special case. We write, for instance,
$C^{(1,2)}_h$
instead of
$C^{(1,i)}_h$,
in accordance with Remark (i) of Section 5.

\subsection*{Computation of the class $z_M$}

Obviously, the class of the cycle
$C^{(i)}_v$
in
$A^1V^{(i)}$
is equal to
$d^{(i)}_vf_i$,
and the class of the cycle
$C^{(i)}_h$
is equal to
$d^{(i)}_hs_i+\beta_if_i$,
where the coefficients satisfy the important inequality
$$
\beta_i\geq d^{(i)}_h
$$
(see Remark
(ii)
in Section 5).
Furthermore, the class of the cycle
$C^{(i,i+1)}_v$
is equal to
$$
d^{(i)}_v(f_i-f_{i+1})
$$
and the class of the cycle
$C^{(i,i+1)}_h$
is equal to
$$
d^{(i)}_hs_i+
\beta_if_i-
(\beta_i-d^{(i)}_h)f_{i+1}.
$$
Setting
$$
\alpha_i=
\left(
(Z^h)^{(i-1)}\cdot L_{i-1}
\right)
$$
in the sense of the definition of the ``intersection index'', which was given
at the
beginning of the paper, we can write down
$$
z^h_i=z^h_{i-1}-\alpha_if_i,
$$
where
$z^h_i$
is the numerical class of the horizontal cycle
$(Z^h)^{(i)}$.

\begin{lll}
The following inequality is true:
$$
\alpha_i\leq\mathop{\rm deg} Z^h=3n^2.
$$
\end{lll}

{\bf Proof.}
Since
$L\subset F$,
and
$\mathop{\rm deg} Z^h$
is equal to
$(Z^h\cdot F)$,
this is obvious. Q.E.D.

\begin{ppp}
The classes
$z_i$
satisfy the following chain of relations:
$$
z_i=z_{i-1}-(2\lambda_in+\lambda^2_i)f_i-\lambda^2_is_i.
$$
\end{ppp}

{\bf Proof.}
We just compute:
$$
z_i=
(D^{(i)})^2=
(D^{(i-1)}-\lambda_iE^{(i)})^2=
$$
$$
=
z_{i-1}-2\lambda_i(D^{(i-1)}\cdot L_{i-1})f_i-
\lambda^2_i(s_i+f_i).
$$
It follows from what was proved in Section 5 that for any
$j\in{\bf Z}_+$
$(D^{(j)}\cdot L_j)=(D\cdot L)=n$.
Q.E.D.

\begin{ppp}
For
$i\geq 2$
in the non-special and for
$i\geq 3$
in the special case the integers
$k_i$,
$\alpha_i$,
$\beta_i$
and
$d^{(i)}_{h,v}$
satisfy the following system of relations:
$$
d^{(i)}_v+\beta_i=
\alpha_i+
d^{(i-1)}_v+
(\beta_{i-1}-
d^{(i-1)}_h)-
2\lambda_in-\lambda^2_i.
$$
For
$i=1$
both in the non-special and special cases we get
$$
d^{(1)}_v+\beta_1=
\alpha_1+
(C^{(0)}_h\cdot L)-
2\lambda_1n-\lambda^2_1,
$$
whereas for
$i=2$
in the special case we get
$$
d^{(2)}_v+\beta_2=
$$
$$
=\alpha_2+d^{(1)}_v+
(\beta_1-
d^{(1)}_h)+
(C^{(0,1)}_h\cdot L_1)-
2\lambda_2n-\lambda^2_2.
$$
\end{ppp}

{\bf Proof.}
To obtain this proposition, it is necessary to write out explicitly the class
of the cycle
$Z^{(i)}$
in terms of the parameters introduced above, and to use the previous
proposition. The
corresponding computations are elementary.

\begin{ppp}
For any
$i\geq 1$
in the non-special case we get the inequality
$$
d^{(i)}_v+\beta_i\leq
(C^{(0)}_h\cdot L)+
\sum^i_{j=1}(3n^2-2\lambda_jn-\lambda^2_j).
$$

In the special case for
$i\geq 2$
we get
$$
d^{(i)}_v+\beta_i\leq
(C^{(0)}_h\cdot L)+
(C^{(0,1)}_h\cdot L_1)+
\sum^i_{j=1}(3n^2-2\lambda_jn-\lambda^2_j),
$$
and
$$
d^{(1)}_v+\beta_1\leq
(C^{(0)}_h\cdot L)+
3n^2-2\lambda_1n-\lambda^2_1.
$$
\end{ppp}

{\bf Proof.}
It is necessary to apply the corresponding inequality of the previous
proposition
$i$
times and to use the last lemma. Q.E.D.

\subsection*{Completing the proof: case A)}

In the case
A)
it is clear that among all the curves, lying on the divisor
$$
\mathop{\bigcup}\limits^M_{i=0}E^{(i,M)},
$$
only those can possibly contain the point
$x_M$,
which lie entirely in
$E^{(M)}$
and are different from the exceptional section
$L_M$.
Consequently, we are justified in writing down
$$
Z^{(M)}=
(Z^h)^{(M)}+
C^{(M)}_v+
C^{(M)}_h+\dots,
$$
where the dots stand for the sum of all the curves, which do not contain the
point
$x_M$.
Set
$$
W=C^{(M)}_v+C^{(M)}_h,
$$
$$
m^v_i=\mathop{\rm mult}\nolimits_{B_{i-1}}W^{i-1},
$$
$$
m^h_i=\mathop{\rm mult}\nolimits_{B_{i-1}}(Z^h)^{(M),i-1}
$$
for
$i\leq L$,
so that
$m_i=m^v_i+m^h_i$.
Obviously, the multiplicities
$m^v_i$
vanish for
$N+1\leq i\leq L$.
Furthermore,
$m^{h,v}_i\leq m^{h,v}_1$,
and similarly to Lemma 6.1 we get
$m^h_1\leq 3n^2$.
Finally,
$m^v_1\leq d^{(M)}_v+d^{(M)}_h\leq d^{(M)}_v+\beta_M$,
so that, summing up our information, we get
$$
3n^2\Sigma_0+
\Sigma_*\left(
(C^{(0)}_h\cdot L)+
\sum^M_{i=1}(3n^2-2\lambda_in-\lambda^2_i)
\right)\geq
$$
$$
\geq
\sum^L_{i=1}p_im^h_i+
\sum^N_{i=1}p_im^v_i\geq
$$
$$
\geq\frac{
\left(2\Sigma_0n+\Sigma_1n+
\varepsilon\sum\limits^M_{i=1}(n-\lambda_i)+
e\right)^2}{
\Sigma_0+\Sigma_1},
$$
if the case is the non-special one. In the special case for
$M\geq 2$
one should add
$(C^{(0,1)}_h\cdot L_1)$
to
$(C^{(0)}_h\cdot L)$.

Now let us consider the non-special case. Replacing
$\Sigma_*$
by
$\varepsilon=\Sigma^*$,
we make our inequality stronger, and replacing
$\varepsilon(C^{(0)}_h\cdot L)$
by
$4ne$,
we get a strict inequality. Subtract the left-hand side from the right-hand one
and look at the
expression just obtained as a quadratic form in
$\lambda_i$
on the domain
$0\leq\lambda_i\leq n$.
By symmetry, its minimum is attained  somewhere on the diagonal line, that is,
at
$\lambda_i=\lambda,$
$0\leq\lambda\leq n$.
Replace all the
$\lambda_i$'s
by this value
$\lambda$.
Thus we get the strict inequality
$$
\Phi<0,
$$
where the expression
$\Phi$
by means of elementary arithmetic can be transformed as follows:
$$
\Phi=(
\Sigma^2_0+
\Sigma_0
\Sigma_1+
\Sigma^2_1)n^2+M\varepsilon
\Sigma_0(n-\lambda)^2-
$$
$$
-M\varepsilon
\Sigma_1(n-\lambda)(n+\lambda)+
$$
$$
+M^2\varepsilon^2(n-\lambda)^2-2e
\Sigma_1n+2M\varepsilon e(n-\lambda)+e^2.
$$

Since
$\lambda\leq n$,
we can replace
$(n+\lambda)$
by
$2n$,
preserving the strict inequality. However, it is easy to check that the last
expression
is the sum of the complete square
$$
(\Sigma_1n-M\varepsilon(n-\lambda)-e)^2
$$
and a few non-negative components. Thus it can not be negative. Our proof is
complete
(in the case under consideration).

In the special case for
$M\geq 2$
the arguments are to be produced in accordance with the same scheme. Just take
into account that here
$\nu(F)=2\varepsilon$,
so that now we may replace the expression
$$
\varepsilon\left((C^{(0)}_h\cdot L)+(C^{(0,1)}_h\cdot L_1)\right)
$$
in the left-hand side by
$4ne$.
The rest part of the computations is the same as in the previous case. For
$M=1$
the computations are much easier.

Q.E.D. for the case A).

\subsection*{Completing the proof: case B)}

As above, we shall trace all the details in the non-special case only. On one
hand, we have the inequality
$$
d^{(M)}_v\leq
(C^{(0)}_h\cdot L)+
\sum^M_{i=1}(3n^2-2\lambda_in-\lambda^2_i).
$$
On the other hand, the following estimate holds:
$$
d^{(M)}_v\geq
\mathop{\rm mult}\nolimits_BZ^{(M)},
$$
whereas for the last multiplicity, in its turn, a lower bound was obtained
above by means of the Iskovskikh-Manin's techniques.
It is easy to see that
$$
\sum^K_{i=1}p^2_i\leq p_1\Sigma\leq\varepsilon\Sigma.
$$
As it was done above, we may assume that all the
$\lambda_i$'s
are equal to
$\lambda$,
$0\leq\lambda\leq n$.
Replacing
$\varepsilon(C^{(0)}_h\cdot L)$
by
$4ne$,
we get the strict inequality
$$
0>
\Sigma^2n^2-2\Sigma
ne+M^2\varepsilon^2(n-\lambda)^2+e^2+2M\varepsilon(n-\lambda)e-
$$
$$
-M\varepsilon(n-\lambda)\Sigma(n+\lambda),
$$
which will be still true when we replace
$(n+\lambda)$
by
$2n$.
But the final expression is a complete square:
$$
(\Sigma n-M\varepsilon(n-\lambda)-e)^2.
$$

This contradiction proves the theorem (in the case under consideration).

In the special case we proceed in the same manner. Here we must add
$(C^{(0,1)}_h\cdot L_1)$
to
$(C^{(0)}_h\cdot L)$.
However, the multiplicity
$\nu(F)=2\varepsilon$
is twice bigger now, so that
eventually we come to the same strict inequality. This contradiction
completes the proof in the case
B).

\subsection*{Completing the proof: case C)}

Let us assume at first that either the case is the non-special one, either it
is special and
$M\geq 3$,
or, finally, that it is special,
$M=2$,
but the point
$x=Z(V,\nu)$
is not the singular point of the fiber
$F$.

Here it is clear that among the curves lying on the divisor
$$
\mathop{\bigcup}\limits^M_{i=0}E^{(i,M)},
$$
only those ones can pass through the point
$x_M$,
which either lie entirely in
$E^{(M)}$
and are different from the exceptional section
$L_M$
(exactly as it was in the case
A)),
or lie entirely in
$E^{(M-1,M)}$
and are different from the exceptional section
$E^{(M-1,M)}\cap E^{(M)}$
(which was already counted in the first group).
Thus
$$
Z^{(M)}=
(Z^h)^{(M)}+W_-+W_++\dots,
$$
where
$W_+$
stands for the sum of all the curves in
$E_+=E^{(M)}$,
which are different from
$L_M$,
$W_-$
stands for the 1-cycle
$C^{(M-1,M)}_v+C^{(M-1,M)}_h$,
and the dots stand for the sum of all the rest curves, which do not pass
through
$x_M$.
Let the symbol
$m^h_i$
mean the same as in the case
A),
and set
$$
m^{\pm}_i=
\mathop{\rm mult}\nolimits_{B_{i-1}}W^{i-1}_{\pm},
$$
$i=1,\dots,L$,
so that
$m_i=m^+_i+m^-_i+m^h_i$.
Obviously, the multiplicities
$m^{\pm}_i$
vanish for
$N_{\pm}+1\leq i\leq L$.
Similarly to the case
A),
we get
$m^{h,\pm}_i\leq m^{h,\pm}_1$,
$m^h_1\leq 3n^2$,
$m^+_1\leq d^{(M)}_v+\beta_M$,
$m^-_1\leq d^{(M-1)}_v+\beta_{M-1}$,
so that finally we come to the following inequality:
$$
3n^2\Sigma_0+
\Sigma_+\left(
(C^{(0)}_h\cdot L)+
\sum^M_{i=1}(3n^2-2\lambda_in-\lambda^2_i)
\right)+
$$
$$
+\Sigma_-\left(
(C^{(0)}_h\cdot L)+
\sum^{M-1}_{i=1}(3n^2-2\lambda_in-\lambda^2_i)
\right)\geq
$$
$$
\geq
\sum^L_{i=1}p_im^h_i+
\sum^{N_+}_{i=1}p_im^+_i+
\sum^{N_-}_{i=1}p_im^-_i\geq
$$
$$
\geq\frac{1}{\Sigma_0+\Sigma_1}
\left(2\Sigma_0n+\Sigma_1n+
\varepsilon_-\sum\limits^{M-1}_{i=1}(n-\lambda_i)+
\varepsilon_+\sum\limits^M_{i=1}(n-\lambda_i)+
e\right)^2,
$$
provided that our case is the non-special one. In the special case for
$M\geq 2$
one should add
$(C^{(0,1)}_h\cdot L_1)$
to
$(C^{(0)}_h\cdot L)$.

Let us consider the non-special case. Replacing
$\Sigma_{\pm}$
by
$\varepsilon_{\pm}=\Sigma^*_{\pm}$,
we preserve the inequality, and replacing
$(\varepsilon_++\varepsilon_-)(C^{(0)}_h\cdot L)$
by
$4ne$,
we make it into a strict one. Subtract the left-hand side from  the right-hand
side and look at the
resulting expression as a quadratic form in the two groups of variables, that
is,
$\lambda^+_i$
and
$\lambda^-_i$,
where we replace
$\lambda_i$
by
$\lambda^{\pm}_i$
in accordance with the following rule: if a variable comes into the sum
$\sum\limits^{M-1}_{i=1}$,
then we replace it by
$\lambda^-_i$,
and if it comes into
$\sum\limits^{M}_{i=1}$,
then we replace it by
$\lambda^+_i$.
The new variables take their values in the domain
$0\leq \lambda^{\pm}_i\leq n$.
By symmetry, the minimum of this quadratic form is attained at some point on
the diagonal plane, that is, at
$\lambda^{\pm}_i=\lambda_{\pm}$,
$0\leq \lambda_{\pm}\leq n$.
Now replace all the
$\lambda^{\pm}_i$'s
by
$\lambda_{\pm}$.
The inequality is still strict. Thus we get
$$
\Phi<0,
$$
where the expression
$\Phi$
can be transformed by means of elementary arithmetic in the following way,
where we set for conveniency
$M-1=M_-$,
$M=M_+$:

$$
\Phi=(
\Sigma^2_0+
\Sigma_0
\Sigma_1+
\Sigma^2_1)n^2+M_-\varepsilon_-
\Sigma_0(n-\lambda_-)^2+
M_+\varepsilon_+
\Sigma_0(n-\lambda_+)^2-
$$
$$
-M_-\varepsilon_-
\Sigma_1(n-\lambda_-)(n+\lambda_-)
-M_+\varepsilon_+
\Sigma_1(n-\lambda_+)(n+\lambda_+)+
$$
$$
+(M_-\varepsilon_-(n-\lambda_-)+
M_+\varepsilon_+(n-\lambda_+))^2-
$$
$$
-2e\Sigma_1n+
$$
$$
+2M_-\varepsilon_-(n-\lambda_-)e
+2M_+\varepsilon_+(n-\lambda_+)e
+e^2.
$$
Since
$\lambda_{\pm}\leq n$,
we can replace
$(n+\lambda_{\pm})$
by
$2n$,
preserving the strict inequality. Now it is easy to check, that the last
expression is the sum of the complete square

$$
(\Sigma_1n-
M_-\varepsilon_-(n-\lambda_-)-
M_+\varepsilon_+(n-\lambda_+)
-e)^2
$$
and a few non-negative components. Again we come to a contradiction, completing
our proof in the case under consideration.

In the special case we use the same arguments, taking into account the equality
$\nu(F)=2(\varepsilon_-+\varepsilon_+)$.

If, finally, our case is the special one with
$M=1$
and
$x=Z(V,\nu)$
is not the singular point of the fiber
$F$,
then the previous arguments work with simplifications.

The only case, which is yet to be considered, is the special one with
$M=1$
or 2, when the point
$x$
is the singularity of the fiber. The case
$M=1$
is more simple. If
$M=2$,
then the point
$x_M$
is the only common point of the following three divisors:
$$
x_2=F^{(2)}\cap E^{(1,2)}\cap E^{(2)}
$$
(the intersection is transversal). Respectively,
$$
Z^{(2)}=
(Z^h)^{(2)}+
C^{(0,2)}_h+
C^{(1,2)}_h+
C^{(1,2)}_v+
C^{(2)}_h+
C^{(2)}_v+
k_2L_2,
$$
where all the 1-cycles but the last one can contain the point
$x_2$.

This case is the only one, when our previous arguments formally do not work
(because of the additional input, which is given by the 1-cycle
$C^{(0,2)}_h\subset F^{(2)}$).
Nevertheless, the general scheme of arguments, which was used in the cases
A)-C)
above, works here, too. We just outline the principal steps of the proof.

Preserving the previous notations, set
$$
N^*=\max\{i|1\leq i\leq K,B_{i-1}\subset E^{i-1}\},
$$
$$
N=\min\{N^*,L\},
$$
$$
\Sigma^*=\varepsilon=\sum\limits^{N^*}_{i=1}p_i,
$$
$$
\Sigma_*=\sum\limits^N_{i=1}p_i\leq\Sigma^*.
$$
Set also
$$
m^0_i=
\mathop{\rm mult}\nolimits_{B_{i-1}}(C^{(0,2)}_h)^{i-1}
$$
for
$1\leq i\leq L$.
Obviously,
$m^0_i=0$
for
$i\geq N+1$.

Now we get the following representation:
$$
\sum^{L}_{i=1}p_im_i=
\sum^{L}_{i=1}p_im^h_i+
\sum^{N_-}_{i=1}p_im^-_i+
\sum^{N_+}_{i=1}p_im^+_i+
\sum^{N}_{i=1}p_im^0_i.
$$
For the four components in the right-hand side we get the following upper
bounds:
$$
\leq 3\Sigma_0n^2,
$$
$$
\leq\varepsilon_-\left(
(C^{(0)}_h\cdot L)+(3n^2-2n\lambda_1-\lambda^2_1)\right),
$$
$$
\leq\varepsilon_+\left(
(C^{(0)}_h\cdot L)+
(C^{(0,1)}_h\cdot L_1)+
\sum^2_{i=1}(3n^2-2n\lambda_i-\lambda^2_i)\right),
$$
$$
\leq\varepsilon m^0_1\leq\varepsilon
(C^{(0,1)}_h\cdot L_1).
$$
It is because of the fourth component that this case is not embraced formally
by the previous arguments. However, here the
multiplicity
$$
\nu(F)=2\varepsilon_++\varepsilon_-+\varepsilon
$$
increases, too. Thus we are able again to replace all the components, into
which
$(C^{(0)}_h\cdot L)$
and
$(C^{(0,1)}_h\cdot L_1)$
come, by
$4ne$.
{}From now on we can just repeat the arguments, which were used in the
``regular'' case C). For
$M=1$
our computations work with considerable simplifications.

The proof of our theorem is complete.
\newpage

\centerline{\large \bf References}

1. Iskovskikh V.A., On the rationality problem for three-dimensional algebraic
varieties, fibered into Del Pezzo
surfaces,-- Proc. of Steklov Math. Inst., V. 208, 1995, 113-122.

2. Iskovskikh V.A. and Manin Yu.I., Three-dimensional quartics and
counterexamples to the L\"uroth problem,--
Math. USSR Sb., V. 86, 1971, 140-166.

3. Iskovskikh V.A. and Pukhlikov A.V., Birational automorphisms of
multi-dimensional algebraic varieties,--
Cont. Math. and Its Appl.,  V. 19, 1995, 3-96.

4. Manin Yu.I., Cubic forms: Algebra, geometry, arithmetic.-- Amsterdam: North
Holland, 1986.

5. Pukhlikov A.V., A remark on the theorem of V.A.Iskovskikh and Yu.I.Manin on
the three-dimensional
quartic,-- Proc. of Steklov Math. Inst., V. 208, 1995, 244-254.

6. Pukhlikov A.V., Birational automorphisms of three-dimensional Del Pezzo
fibrations,-- Warwick Preprint 30/1996.

7. Pukhlikov A.V., Essentials of the method of maximal
singulariti\-es,-- Warwick Preprint 31/1996.

8. Sarkisov V.G., Birational automorphisms of conic bundles,-- Math. USSR
Izvestia, V. 17, 1981, No. 1, 177-202.

9. Sarkisov V.G., On conic bundles structures,-- Math. USSR Izvestia, V. 20,
1982, No. 2, 355-390.

\end{document}